\documentclass{amsart}

\usepackage{amsmath,amssymb}
\usepackage{graphicx}
\usepackage{enumerate}

\makeatletter
\@addtoreset{equation}{section}
\makeatother

\numberwithin{equation}{section}

\def\ii{{\sqrt{-1}}}

\def\SO{\mathrm {SO}}

\def\Ind{\mathrm {Ind}}
\def\Res{\mathrm {Res}}

\def\fri{\mathfrak {i}}
\def\fj{\mathfrak {j}}
\def\fm{\mathfrak {m}}

\def\FF{\mathcal {F}}
\def\LA{\langle}
\def\RA{\rangle}

\def\EE{{\mathbb{E}}}
\def\CC{{\mathbb C}}
\def\RR{{\mathbb R}}

\def \even{\mathrm {even}}

\def \dalpha{{\dot \alpha}}
\def \dbeta{{\dot \beta}}

\def \Cliff{{\mathrm{Cliff}}}
\def \CLIFF{{\rm{CLIFF}}}
\def \Gr{{\rm{Gr}}}
\def \CLIFFC{{\mathrm{CLIFF}^{\mathbb C}}}
\def \CG{{\mathrm{CG}}}

\def \END{{\mathrm{END}}}

\def \Hom{{\mathrm{Hom}}}
\def \Cinf{{{\mathcal C}^\infty}}

\def \Ker{{\mathrm{Ker}}}

\def \ii{{\sqrt{-1}}}
\def\Not{\not \!\!}

\def\pt{{pt}}
\def\iotaS{{\iota_{M, S}}}
\def\iotaRS{{\iota_{\RR^n, S}}}
\def\iotaRM{{\mu_{\RR^n, M}}}
\def\iotaMTS{{\iota_{M, T_S}}}
\def\piTS{{\pi_{S, T_S}}}
\def\overwave{\tilde}

\def \sa{{\mathrm{sa}}}

\newtheorem{definition}{Definition}[section]

\newtheorem{theorem}{Theorem}[section]
\newtheorem{proposition}{Proposition}[section]

\newtheorem{remark}{Remark}[section]
\newtheorem{lemma}{Lemma}[section]

\def\book#1{\rm{#1}, }
\def\paper#1{\textit{#1}, }
\def\jour#1{\rm{#1}, }
\def\yr#1{({\rm{#1}) }}
\def\vol#1{\textbf{#1}}
\def\pages#1{\rm{#1}}

\def\publ#1{\rm{#1}, }
\def\by#1{{\rm{#1}, }}

\begin{document}

{\large{
\title{
Generalized Weierstrass Relations
and Frobenius reciprocity}

\author{Shigeki Matsutani}

\date{Nov 9, 2006}

\address{
Shigeki Matsutani\newline
}


\maketitle

\begin{abstract}
This article investigates local properties of
the further generalized Weierstrass relations
for a spin manifold $S$ immersed in a higher dimensional spin
manifold $M$ from viewpoint of study of submanifold quantum
mechanics. We show that kernel of a certain 
Dirac operator defined over $S$, which
we call submanifold Dirac operator, gives the data of
the immersion. In the derivation, 
 the simple Frobenius reciprocity of Clifford algebras $S$
and $M$ plays important roles.
\end{abstract}



\section{ Introduction}

This article is a sequel of the previous paper \cite{Ma4}.
We study a connection between
the generalized Weierstrass relation
and Frobenius reciprocity, which is partially described
in \cite{Ma4}, and obtain
further generalized Weierstrass relation
over a spin manifold $S$ immersed in higher dimensional
spin manifold $M$.

The generalized Weierstrass relation is a generalization of the
Weierstrass relation appearing in the minimal
surface theory \cite{E}, which gives data of
immersion of a conformal surface in higher dimensional
flat spaces, e.g., euclidean space.
Although similar relations appeared in \cite{E} and it was
obtained by K. Kenmotsu \cite{Ke},
the generalized Weierstrass relation was mainly studied in 1990's,
by B. G. Konopelchenko \cite{Ko} and G. Landolfi \cite{KL},
F. Pedit and U. Pinkall \cite{PP}, I. A. Taimanov \cite{Ta1}, 
and so on.
Their studies are, basically, in the framework of geometrical 
interpretations of integrable system. 
In the studies a certain Dirac operator appears and its global
solutions of its Dirac equation provides the data of immersion of surfaces;
in this article, we shall call, later, the Dirac operator (equation)
{\it{submanifold Dirac operator (equation)}}.
T. Friedrich investigated the relations for a surface
immersed in the euclidean 3-space $\RR^3$ from a viewpoint of 
study of Dirac operator
\cite{Fr}. V. V. Varlamov also studied the relations from a point of
view of Clifford algebra \cite{V}.
In \cite{KL}, surfaces in flat $n$-space are treated and
those in Riemann spaces were mentioned. Further L. V. Bogdanov and
E. V. Ferapontov generalized the relation to a surface in projective
space \cite{BF}. Recently
I. A. Taimanov gave a proper survey on the related topics 
and open problems in \cite{Ta3}. 

The author has studied the submanifold Dirac operator
since 1990 in the framework of 
quantum mechanics over submanifolds,
which we call submanifold quantum mechanics;
in the framework, we deal with a restriction of differential
operator, hamiltonian, defined over a manifold to one over its submanifold
and then we find a non-trivial structure in the operator due to
the half-density \cite{Ho, Ma5}.
In \cite{MT, Ma1} he and his coauthor investigated the Dirac operator over  curves
in flat space and 
showed that the Dirac operator is identified with the
 operator of the Frenet-Serret relation and 
a natural linear operator in the soliton theory.
The latter one gives a geometrical interpretation of integrable system.
When we apply the scheme developed in \cite{MT, Ma1} to
the immersed surface case [\cite{Ma2} and reference therein],
we also encounter the same situation;
the Dirac operator coincides with the Dirac operator appearing in
the generalized Weierstrass relations and with a natural linear operator 
of an two-dimensional soliton equation.
Further the analytic torsion of the submanifold Dirac operator 
is also connected with globally geometrical properties  \cite{Ma2, Ma3},
as the Dirac operator with gauge fields exhibits geometrical
properties of its related principal bundle via the analytic torsion
in the framework of the Atiyah-Singer index theorem and so on \cite{BGV};
the submanifold Dirac operator is also directly associated with the
global geometry.

In the series of works, the author has considered 
why the Dirac operator given in the framework of
the submanifold quantum mechanics
appears in the generalized Weierstrass relation
 and expresses geometrical properties of submanifold.
In other words, our motivation of the study is to clear
what is the submanifold quantum mechanics and what is 
the generalized Weierstrass relation from viewpoint of
study of the submanifold quantum mechanics.

In fact, recently shape effect in quantum mechanics becomes
to play a more important role in physics due to development of nanotechnology.
The submanifold Schr\"odinger operator in the submanifold quantum
 mechanics is applied to more actual geometrical objects 
\cite{En, MBK, GWL, MEE}.
Thus it is required to reveal mathematical (analytic, geometrical
and algebraic) structure of the submanifold quantum mechanics.
On the submanifold Schr\"odinger operator,
its algebraic essential was clarified \cite{Ma5}.

This article is the final version of
the studies on the construction and the local properties of 
the submanifold Dirac operators. 
We find the answers to the problem why the submanifold Dirac operator 
constructed in the framework of the submanifold quantum mechanics
represents immersed geometry; 
this means a local aspect of
the generalized Weierstrass relation.
Though, of course, the global feature of the Dirac equation might
be more interesting than local ones,
its essential of the answer is based on local properties,
which are connected with the simple Frobenius relations in 
a local chart.
Thus it is not difficult to generalize the submanifold
Dirac operator defined over a surface immersed in $\RR^n$ to
one over more general geometrical situations, at least locally.
If there is no obstruction, it might determine a submanifold globally.

Here we note that the definition and construction of
our submanifold Dirac operator differs from that of C. B\"ar \cite{B},
though both forms coincide. On the line, 
N. Ginoux and B. Morel \cite{GM}, and H. d'Oussama and X. Zhang \cite{OZ}
also investigated eigenvalues of the submanifold Dirac operator.
However their construction is not directly 
associated with our requirement.
Thus we concentrate ourselves into the reveal using the
scheme of the submanifold quantum mechanics.

We will mention our plan of this article. Section 2 shows our
conventions of the Clifford algebra whereas section 3 provides
our geometrical assumptions and conventions of this article.
After we consider the Dirac operator over a manifold in section 4,
we will construct a Dirac operator over its submanifold and investigate it
in section 5. There we will give our main theorem as
Theorem 5.1.

\bigskip

\section{Local expression of Clifford Algebra}

In this section, in order to show our convention in this article,
we will briefly review  the Clifford algebra
\cite{ABS, C, GW}.
The Clifford Algebra $\CLIFF(\RR^m)$ is introduced as
a quotient ring of a tensor algebra, $\mathbb{T}(\RR^m)/((v,u)_{\RR^m}-1)$,
where
$u$, $v$ are elements of $m$-dimensional vector space
$\RR^m$ and $(v,u)_{\RR^m}$ is the natural inner product.

With respect to the degree of a tensor product, we have a natural
filtration $\FF^{\ell} \CLIFF(\RR^m) \supset \FF^{\ell-1} \CLIFF(\RR^m)$,
where $\FF^{0} \CLIFF(\RR^m) = \RR$
 and $\FF^p\CLIFF(\RR^m) = 0$ for $p <0$,
with a graded algebra
$\CLIFF^p(\RR^m):=\FF^p \CLIFF^1(\RR^m)/ \FF^{p-1} \CLIFF^1(\RR^m)$.
Let its subalgebra with even degrees be denoted by
$\CLIFF^{\even}(\RR^m) = \cup_{p=\even}^m \CLIFF^p(\RR^m)$.

The exterior algebra $\wedge \RR^m= \oplus_{j=1}^m \wedge^j \RR^m$,
is isomorphic to $\CLIFF(\RR^m)$ as $\RR^m$ vector space,
$\wedge^p \RR^m \to \CLIFF^p(\RR^m)$
and thus let the isomorphism,
\begin{gather}
 \gamma^{(m)}: \RR^m \to \CLIFF^1(\RR^m).
\label{eq:gammaRm}
\end{gather}
For the basis of $\RR^m$ denoted by $( e^{(m), i})_{i=1, \cdots, m}$,
let $*$ operator be the involution in $\CLIFF(\RR^m)$ such that
$(\gamma^{(m)}(e^{(m), i_1}) $ $\cdots \gamma^{(m)}(e^{(m), i_j}))^*:=$
$(\gamma^{(m)}(e^{(m), i_j}) $ $\cdots \gamma^{(m)}(e^{(m), i_1}))$.

Let  $\Cliff(\RR^m)$ be a left $\CLIFF(\RR^m)$-module whose
endomorphism \break $\END(\Cliff(\RR^m))$ is isomorphic to $\CLIFFC(\RR^m)$
$(\equiv\CLIFF(\RR^m)\otimes \CC)$ as
$2^{[n/2]}$ dimensional $\CC$-vector space representation;
$\epsilon_m : \CLIFFC(\RR^m) \to \END(\Cliff(\RR^m))$.
Let  $\Cliff^*(\RR^m)$ be a right $\CLIFF(\RR^m)$-module which
is isomorphic to $\Cliff(\RR^m)$;
$\varphi: \Cliff(\RR^m) \to \Cliff(\RR^m)^*$;
for $C\in \CLIFF(\RR^m)$ and $c\in \Cliff(\RR^m)$,
$\varphi(C c) = \varphi(c) C^*$ and let $\overline{c}:=\varphi(c)$.

We may find bases $(c^{(m), a})_{a=1, \cdots, 2^{[m/2]}}
\in \Cliff(\RR^m)$ such that
for $\overline{c^{(m), a}} $ $= \varphi(c^{(m), a})$,
 $\overline{c^{(m), a}} c^{(m), b} = \delta_{a,b}$.
Every $\psi^{(m)} \in \Cliff(\RR^m)$ is expressed as
$\psi^{(m)}=\sum_{a=1}^{2^{[m/2]}} \psi_a^{(m)} c^{(m), a}$.
For $\overline{\phi^{(m)}}
=\sum_{a=1}^{2^{[m/2]}} \overline{\phi_a^{(m)}} \ \overline{c^{(m), a}}$
and  ${\psi^{(m)}}$,
we will introduce a natural pairing:
\begin{gather}
\LA  \RA_{\Cliff(\RR^m)}: 
 \Cliff(\RR^m)^*\times \Cliff(\RR^m) \to \CC,
\label{eq:pairing}
\end{gather}
by
$$
\LA  \overline{\phi^{(m)}}, \psi^{(m)}\RA_{\Cliff(\RR^m)}=
\sum_{a=1}^{2^{[m/2]}} \overline{\phi_a^{(m)}} \psi_a^{(m)}.
$$
For multiplicative group of $\CLIFF^{\even}(\RR^m)$,
$\CLIFF^{\even,\times}(\RR^m)$,
the Clifford group $\CG(\RR^m)$ is defined by
$$
	\{ \tau \in \CLIFF^{\\even\times}(\RR^m)\ |
         \ \text{ for }\forall v \in \CLIFF^1(\RR^m),
          \tau v\tau^* \in \CLIFF^1(\RR^m) \}.
$$
For representations $\epsilon_m$ and $\epsilon_m'$, there exists
$\tau \in \CG(\RR^m)$ and an action $A_\tau$ on $\epsilon$'s such that
$A_\tau \epsilon_m(C) = \epsilon_m'(\tau C \tau^{-1})$
for $C \in \CLIFF(\RR^m)$.
Due to $\gamma^{(m)} : \RR^m \to \CLIFF^1(\RR^m)$ and 
(\ref{eq:pairing}), we have
\begin{equation}
\LA\ ,\  \RA_{ \Cliff(\RR^m)}:
\Cliff(\RR^m)^*\times\gamma^{(m)}(\RR^m)\times\Cliff(\RR^m) \to \CC.
\label{eq:2-3}
\end{equation}
This is a linear map from $\RR^m$ to $\CC$.
Let the coproduct be $\mathfrak{m}: 
 \Cliff(\RR^m)\to \Cliff(\RR^m)\times \Cliff(\RR^m)$,
$(\psi \mapsto (\psi, \psi))$.
Restricted domain to its inverse image of $\RR\subset \CC$,
(\ref{eq:2-3}) with the operator
{\lq\lq}$\LA, \gamma^{(m)}(\cdot) \ \RA\circ \varphi \otimes 1\circ 
\fm${\rq\rq} can be regarded as $\Hom_\RR(\RR^m, \RR) \approx \RR^m$.
Hence we have the following lemma.

\begin{lemma} \label{lemma:SOm}
There exists a subset $\Cliff^{pr}(\RR^m)$ of  $\Cliff(\RR^m)$
which is isomorphic to $\RR^m$ as $\RR$-vector space such that
for
$$
\fri: \RR^m \to \Cliff^{pr}(\RR^m) \subset \Cliff(\RR^m),
\quad
(v \mapsto \psi^{pr}_v),
$$
$$
\fj :=\LA\ , \gamma^{(m)}(\cdot)\ \RA \circ \varphi\otimes 1 \circ 
       \fm \otimes 1\circ \fri \otimes 1: \RR^m \times \RR^m \to \RR
$$
is identified with  the inner product $(,)_{\RR^m}: \RR^m \times \RR^m \to \RR$,
$(u, v)_{\RR^m} \equiv \fj(u, v)$, i.e.,
$$
\LA \overline{\psi^{pr}_v}, \gamma^{(m)}(w) \psi^{pr}_v\RA_{\Cliff(\RR^m)}
 = (v, w)_{\RR^m}.
$$
\end{lemma}

This lemma shows that
there exist elements $( \phi_{e^{(m),j}}^{(m)})_{j=1, \cdots, m}$ of
 $\Cliff^{pr}(\RR^m)$;
 for $b^{(m), i} =\sum_{j=1}^m \Lambda^i_{\ j} e^{(m), j}$,
\begin{gather}
\LA\overline \phi_{e^{(m), \ell}}^{(m)}, \gamma^{(m)}(b^{(m), i})
 \phi_{e^{(m),\ell}}^{(m)}
\RA_{\Cliff(\RR^m)} = \Lambda^i_{\ \ell}.
\label{eq:RmSOm}
\end{gather}
This correspondence is well-known in physicists,
which is, of course, independent from
the coordinate system and gives the data of $\SO(\RR^m)\times\RR$.

\bigskip

Here let us consider an
embedding $\RR^k$ into $\RR^n$  ($k<n$):
$\iota_{n, k}:\RR^k \hookrightarrow \RR^n$
and
$\pi_{k, n}:\RR^n \longrightarrow \RR^k$ such that
for $u^{(n)} \in \RR^n$ and $v^{(k)} \in \RR^k$,
$(\iota_{n, k} u^{(n)}, v^{(k)})_{\RR^{n}}$ 
$\equiv (u^{(n)}, \pi_{n,, k}v^{(k)})_{\RR^{k}}$.
In this article, we are concerned with the 
moduli of the embedding or Grassmann manifold
$\Gr_{n, k}:= \SO(n)/\SO(k) \SO(n-k)$.
The embedding $\iota_{n, k}$ corresponds to a point $q$ of
$\Gr_{n, k}$  $= \SO(\RR^n)/ \SO(\RR^n)  \SO(\RR^{n-k})$.
Using the Clifford module we will deal with them like \cite{ABS}.
The following proposition is obvious due to \cite{ABS}.
\begin{proposition} \label{prop:CLkn}
\begin{enumerate}
\item
For $k < n$, $\CLIFF(\RR^k)$ is a subalgebra $\CLIFF(\RR^n)$
by the natural inclusion of generators.
$\iota_{n, k}^{\flat}: \CLIFF(\RR^k) \to \CLIFF(\RR^n)$.

\item
For $k < n$, $\CG(\RR^k)$ is a natural subgroup of $\CG(\RR^n)$.
\end{enumerate}
\end{proposition}

The $\iota_{n, k}$ and
 $\pi_{k, n}$
give an induced representation and
 a restriction representation:
There exists an element $\tau_{q}$ in $\CG(\RR^n)$
such that
\begin{gather*}
\begin{array}{lll}
{\Ind^{\tau_q}}_k^n &: \Cliff(\RR^k) \to \Cliff(\RR^n), &
\left(
\displaystyle{ {\Ind^{\tau_q}}_k^n \psi^{(k)}:=
\sum_{a=1}^{2^{[k/2]}}
 \psi_a^{(k)} \left[\sum_{b=1}^{2^{[n/2]}} {\tau_q}^a_{\ b} c^{(n), b}\right]
}  \right), \\
{\Res^{\tau_q}}_k^n &: \Cliff(\RR^n) \to \Cliff(\RR^k), &
\left(
\displaystyle{ {\Res^{\tau_q}}_k^n \psi^{(n)}:=\sum_{a=1}^{2^{[k/2]}}
 \left[\sum_{b=1}^{2^{[n/2]}}
  \psi_b^{(n)}
  {\tau_q^{-1}}^{\ b}_{a}\right] c^{(k), a}}
\right) .\\
\end{array}
\end{gather*}
The Frobenius reciprocity gives
for  $\psi^{(k)} \in \Cliff(\RR^k)$ and
 $\phi^{(n)} \in \Cliff(\RR^n)$,
\begin{gather}
\LA \overline{\Res^{\tau_q}}_k^n \psi^{(n)}, {\phi^{(k)}} \RA_{\Cliff(\RR^k)}=
\LA \psi^{(n)}, \overline{{\Ind^{\tau_q}}_k^n \phi^{(k)}} \RA_{\Cliff(\RR^n)}.
\label{eq:FrobRec}
\end{gather}

\bigskip

Using the relation (\ref{eq:FrobRec}),
we will consider $\Lambda$ in  (\ref{eq:RmSOm})
and its relation to the point $q$ of the Grassmannian $\Gr_{n, k}$.
For $u^{(n)} \in \RR^n$ and $v^{(k)} \in \RR^k$, 
let $\psi_{u^{(n)}}^{(n)} = \fri(u^{(n)})$
and  $\psi_{u^{(n)}}^{(k)}:={\Res^{\tau_q}}_k^n \psi_{u^{(n)}}^{(n)}$
using  $\tau_q \in \CG(\RR^n)$,
and then we have the relation,
$\gamma^{(n)}(\iota_{n, k} (v^{(k)}))  \psi_{u^{(n)}}^{(n)}
 ={\Ind^{\tau_q}}_k^n \gamma^{(k)}(v^{(k)}) \psi_{u^{(n)}}^{(k)}$.
The Frobenius reciprocity (\ref{eq:FrobRec}) gives
\begin{gather}
\begin{split}
\LA
 \overline{\psi_{u^{(n)}}^{(k)}},\
 \gamma^{(k)}(v^{(k)}) \psi_{u^{(n)}}^{(k)}
\RA_{\Cliff(\RR^k)} 
& =
\LA
 \overline{\psi_{u^{(n)}}^{(n)}},\
 {\Ind^{\tau_q}}_k^n  \gamma^{(k)}(v^{(k)}) \psi_{u^{(n)}}^{(k)}
\RA_{\Cliff(\RR^n)}\\
& =
\LA
 \overline{\psi_{u^{(n)}}^{(n)}},\
    \gamma^{(n)}(\iota_{n, k} (v^{(k)}))  \psi_{u^{(n)}}^{(n)}
\RA_{\Cliff(\RR^n)}\\
& = (\iota_{n, k}(v^{(k)}), u^{(n)})_{\RR^n}.
\label{eq:EK}
\end{split}
\end{gather}
Every pair $(u^{(n)}, v^{(k)})$ recovers the point $q$ in $\Gr_{n, k}$.

This relation (\ref{eq:EK}) has an alternative expression 
using another reference
embedding $\iota_{n, k}^o:\RR^k \hookrightarrow \RR^n$
 associated to a base point of $o \in \Gr_{n, k}$ and $\tau_o \in \CG(\RR^n)$. 
For given $\tau_q$ and $\tau_o$ of $\CG(\RR^n)$,
we find an element $\tau \in \CG(\RR^n)$ such that $\tau_q = \tau^{-1} \tau_o$.
When one wishes to consider $\tau_q$ as a representation of $\Gr_{n,k}$, 
he could deal with the element $\tau$ 
by fixing $\tau_o$;
and investigate $\Gr_{n, k}$.
Then we have
$$
 {\Ind^{\tau_q}}_k^n  = \tau^{-1} {\Ind^{\tau_o}}_k^n,
\quad
 {\Res^{\tau_q}}_k^n  =  {\Res^{\tau_o}}_k^n \tau.
$$
For the situations of (\ref{eq:EK}),
let $\phi_{u^{(n)}}^{(n)}:=\tau \psi_{u^{(n)}}^{(n)}$
and then we have 
$\psi_{u^{(n)}}^{(k)}={\Res^{\tau_o}}_k^n \phi_{u^{(n)}}^{(n)}$.
(\ref{eq:EK}) becomes
\begin{gather}
\begin{split}
\LA
 \overline{\psi_{u^{(n)}}^{(k)}},\
 \gamma^{(k)}(v^{(k)}) \psi_{u^{(n)}}^{(k)}
\RA_{\Cliff(\RR^k)} 
& =
\LA
 \overline{\phi_{u^{(n)}}^{(n)}},\
    \gamma^{(n)}(\iota_{n, k}^o (v^{(k)}))  \phi_{u^{(n)}}^{(n)}
\RA_{\Cliff(\RR^n)}\\
& = (\iota_{n, k}(v^{(k)}), u^{(n)})_{\RR^n}.
\label{eq:ek}
\end{split}
\end{gather}
This also  provides the data of $\Gr_{n, k}$ and the immersion,
which essentially comes from (\ref{eq:FrobRec}) and Lemma 2.1.
We will use latter relation (\ref{eq:ek}) for the generalized
Weierstrass relations.

\section{Geometrical Preliminary}

In this section,
 we will give a geometrical preliminary.
As we use primitive facts in sheaf theory \cite{I},
first we show our conventions  as follows.
For a fiber bundle $A$ over a paracompact
differential manifold $X$ and an
open
set $U \subset X$, let $A_X$ denote a
sheaf given by  a set of smooth local
sections of the fiber bundle $A$, {\it e.g.},
$\CC_X^r$ is a  sheaf given by
smooth local sections of
complex vector bundle over $X$ of rank $r$,
and $A_X(U)\equiv\Gamma(U, A_X)$
sections of $A_X$ over $U$.

Further for open sets $U\subset V\subset X$,
the restriction of a sheaf $A_X$ is denoted 
by $r_{U V}$. 
Using the direct limit
for $\{ U\ |\ \pt\in U \subset X\}$,
we have a {\it stalk} $A_{\pt}$ of $A_X$ by setting
$A_{\pt}\equiv \Gamma(\pt, A_X)
:=\lim_{U \rightarrow \pt}  A_X(U)$.
Similarly for a compact subset $K$ in $X$,
$i_K: K \hookrightarrow X$ and
 for $\{ U\ |\ K\subset U \subset X\}$,
we have
$\Gamma(K, A_X):=\lim_{U \rightarrow K} A_X(U)$ and
  $r_{K, U} A_X$.

On the other hand, for a topological subset $Y$ of
$X$, $i_Y:Y\hookrightarrow X$, there is an inverse sheaf,
$i_Y^{-1} A_X$ given by the sections
$i_Y^{-1} A(U)=\Gamma(i_Y(U), A_X)$ for $U\subset Y$.
When $Y$ is a compact set, 
we have an equality
$\Gamma(i_Y Y, A_X)= \Gamma(Y, i_Y^{-1}A_X)$
 (Theorem 2.2 in \cite{I}) and
 we identify them in this article.
Further $\Gamma_c(U, A_X)$ denotes the set of
smooth sections of $A_X$ whose support is
 compact in $U$. 
For a compact subset $K$ of $X$,
$\Gamma_K(X, A_{X})$
 is a set of global sections
of $A_{X}$
whose support is in $K$.

\bigskip
Let $(M, g_M)$ be  a $n$ spin manifold,
which is acted by a Lie transformation group $G$ as
its isometory.
The metric $g_M$ of $M$ is a global section of sheaf
$\Hom_\RR(\Theta_M, \Omega_M)$,
where $\Theta_M$ and $\Omega_M$  are
sheaves of tangent and cotangent spaces as $\Cinf$-modules:
$g_M(\ , \ ) : \Theta_M \times \Theta_M \to \RR_M$.

Let us consider a locally closed $k$ spin manifold
$S$ embedded in $M$ \cite{W, He}; $\iotaS: S \hookrightarrow M$,  so
that for every point $p$ in $S$,
there is a subgroup $H$ of $G$ satisfying
\begin{equation}
    T_p M = T_p (H \circ p) \oplus T_p S.
\label{eq:TM}
\end{equation}
We  identify $\iotaS(S)$ with $S$.
$H$ may depend on the position $p$ in general.

Since $\iotaS^{-1}\Theta_M$ can be regarded as
a subsheaf of the $(n, k)$ Grassmannian sheaf $\Gr^{(n, k)}_S$
over $S$,  fixing a section $\Gr^{(n, k)}_S$ 
corresponds to determine the immersion $\iotaS$ up to
global symmetry like euclidean moves.
We consider $\iotaS^{-1} \Theta_M$ and $\iotaS_* \Theta_S$.
Let $\Theta_S^\perp := \iotaS^{-1} \Theta_{M} /\Theta_S$;
$\Gr^{(n, k)}_S$ can be realized as
the quotient of orthogonal group sheaves
$\Gr^{(n, k)}_S = \SO(\iotaS^{-1}\Theta_M)/ \SO(\Theta_S) \SO(\Theta_S^\perp)$.

For example, as $r_{S, M}$ is defined by
a direct limit of open sets of $M$ to $\iotaS(S)$.
we should consider its vicinity in $M$.
We prepare a tubular neighborhood $T_S$ of $S$ in $M$;
 $\piTS:T_S \to S$ and   $\iotaMTS: T_S \hookrightarrow M$.

As our theory is local and we use only germs at a point in vicinity of
$\iotaS(S)$,  we consider a sufficiently
small open set $U$ in $M$ such that
 $U\cap S \neq \emptyset$ instead of
$M$ and $S$;
Without loss of generality,  we assume that
$M$ and $S$ are diffeomorphic to $\RR^n$ and $\RR^k$ respectively,
 there exists a compact subset $K$ of $M$ such that $S \subset K$,
 and later we may sometimes identify $M$ with $T_S$.

Further due to the group action $H$, 
we assume that $T_S$ and $S$ satisfy the following conditions.

\begin{enumerate}
\item $T_S$ behaves as a normal bundle $\piTS: T_S \to S$,

\item there exist the base $b^{(n),\dalpha}$ ($\dalpha = k+1, \cdots, n$) 
of $T_S$ and $\Theta_S^\perp$, its dual base ${b^{(n)}_\dalpha}$,
 and $q:=(q^\dalpha)_{\dalpha = k+1, \cdots, n}$ 
 the normal coordinate of $T_S$
such that
1) for $X \in \Theta_S(S)$ and the Riemannian connection $\nabla_X$ in $M$,
 $\nabla_X {b^{(n)}_\dalpha}$ belongs to $\Theta_S(S)$
(See proof of Lemma \ref{lemma:g_S})
and
2) every point $pt \in T_S$
is expressed by $pt =\piTS \pt+q^\dalpha {b^{(n)}_\dalpha}$.

\item $T_S$ and $S$ have local parameterization.
$u: T_S \to \RR^{k}\times \RR^{n-k}$ such that $u=(s, q)$ 
and $s: S \to \RR^k$;
$u=(u^\mu)_{\mu=1, \cdots, n}$
$=(s^\alpha,q^\dalpha)_{\alpha = 1, \cdots, k, \dalpha = k+1, \cdots, n}$.
(We will use the Einstein convention.)
\end{enumerate}

Let $S_q$ be $u^{-1}(\RR^k \times \{q\})$ for fixing $q$.
$\{S_q\}_q$ has a foliation structure.
As a result of (\ref{eq:TM}),  $S$ could be interpreted as
an analytic manifold;
$$
	S\equiv \{(s, q) \in T_S \ |\ q = 0\}.
$$
For every sheaf $A_{T_S}$ of $T_S$,
we have a sheaf $A_S$ of $S$ and a restriction map
$r_{S, M} : A_{T_S}^n \to A_S^n$
by substituting $q=0$ into $f(s, q)$.
Hereafter we use the symbol $r_{S, M}$ in this meaning.
Due to the above assumption,
the metric $g_{T_S}$ of $T_S$ at $(s, q)$ induced from $M$
is given as
\begin{equation}
 g_{T_S}=\begin{pmatrix} g_{ {S_q}} & 0 \\ 0 & 1 \end{pmatrix},
  \label{eq:gTS}
\end{equation}
where $g_{S_q}$ is a metric $S_q$ given by proof 
of Lemma \ref{lemma:g_S}.
We also introduce objects and maps for $S_q$ as for $S$, e.g., $\iota_{M, S_q}$.

\begin{lemma} \label{lemma:g_S}
Let $g_{T_S}$ and $g_S$ be induced metrics of $g_M$ and
$\Gamma_\dalpha/k$ be the mean curvature vector field along
$b^{(n) \dalpha}$
{\rm{[\cite{W} p.119]}},
\begin{equation}
\det g_{T_S}=\rho \det g_S , \quad
\rho= (1 + \Gamma_{\dalpha} q^\dalpha
+ \mathcal O(q^\dalpha q^\dbeta) )^2 . 
\label{2-3}
\end{equation}
\end{lemma}

\begin{proof}
In general, we consider more general normal 
unit vectors $\tilde b^{(n)}_\dalpha$
$\in T_{\pt}^\perp S$ at $\pt \in S$.
At a point in $S$,
we find the Christoffel symbol
$\Gamma^\beta_{\ \dbeta\alpha}$ over $S$
as a relation, the equation of Weingarten [\cite{W} p.119],
\begin{equation}
    \nabla_\alpha \tilde b^{(n)}_\dbeta =
    \Gamma^\beta_{\ \dbeta\alpha}  b^{(n)}_{\beta}
  +\tilde \Gamma^\dalpha_{\ \alpha\dbeta}
  \tilde b^{(n)}_{\dalpha}. 
\label{eq:frenet1}
\end{equation}
Here $\nabla_\alpha$ is the Riemannian connection of $M$ for the direction
$\partial/\partial s^\alpha$ of $TS$, and
$ b^{(n)}_{\beta} := \partial/\partial s^\beta$ of $TS$.
Let $\Lambda_\dalpha^{\ \dbeta}$ be a section of $\SO(\Theta_S^\perp)$
 such that 
its Lie algebraic parameter $\theta_{\dalpha,\dbeta}$ satisfies
$\partial_\alpha \theta_{\dalpha,\dbeta}
=\tilde\Gamma^\dalpha_{\ \alpha\dbeta}$ noting
$\tilde\Gamma^\dalpha_{\ \alpha\dbeta} =-\tilde\Gamma^\dbeta_{\ \alpha\dalpha}$.
Since $S$ is homeomorphic to $\RR^k$, we can find such a parameter 
$\theta_{\dalpha,\dbeta}$ by solving the differential equation.

Let $b^{(n)}_\dalpha = \Lambda_\dalpha^{\ \dbeta} \tilde b^{(n)}_\dbeta$.
Then (\ref{eq:frenet1}) is reduced to
\begin{equation*}
    \nabla_\alpha b^{(n)}_\dbeta =
    \Gamma^\beta_{\ \dbeta\alpha} b^{(n)}_{\beta}.
\end{equation*}

For a point $pt$ in $T_S$, the
 moving frame $e^{(n), i} = d x^i \in \Gamma(pt, \Theta_{T_S})$
is expressed by
$e^{(n),i}_{\ \alpha} d s^\alpha  = (\piTS_*(e^i_{\alpha}) +
            q^\dalpha\Gamma^\beta_{\
            \dalpha\alpha} b^{(n),i}_{\ \beta})d s^\alpha $.
The metric in $T_S$  and its determinant are
given by
\begin{equation}
g_{ {S_q} \alpha\beta}
  = g_{S \alpha\beta}+
    [\Gamma_{\ \dalpha\alpha}^\Gamma g_{S\gamma\beta}+
    g_{S\alpha\gamma}\Gamma_{\ \dalpha\beta}^\gamma]q^\dalpha
    +[\Gamma_{\ \dalpha\alpha}^\delta g_{S\delta\gamma}
     \Gamma_{\ \dbeta\beta}^\gamma]q^\dalpha q^\dbeta ,
        \label{2-2}
\end{equation}
where $g_{S \alpha\beta}:=g_{M\ i,j}e^i_\alpha e^j_\beta$;
Let $\Gamma_{\ \dbeta}:=\Gamma^\alpha_{\ \dbeta\alpha}$ over $S$;
$(\Gamma_{\ \dbeta})/k$ is the mean curvature vector of 
$b^{(n) \dbeta}$ [\cite{W} p.119].
\end{proof}

\bigskip
There is an action of $\SO(\Theta_S^\perp)$ on $\Theta_S^\perp$.
Obviously (\ref{2-3}) is invariant for the action
 $\SO(\Theta_S^\perp)$.

\bigskip

\section{Dirac System in $M$}

For the above geometrical situation,
we will consider a Dirac equation over $M$ [\cite{BGV} 3.3] here.

We, first, introduce a paring given by the pointwise product
$\LA,\RA_{\Cliff_M}$
for the germs of the Clifford module $\Cliff_{M}$  over $M$ and
its natural hermite conjugate $\Cliff_{M}^*$;
$\varphi_{\pt}$ is the  hermite conjugate operator which gives
the isomorphism from  $\Cliff_{M}$ to $\Cliff_{M}^*$ and
$\LA \overline\psi_{M, 1} \psi_{M, 2}\RA_{\Cliff_M}
 \in \Gamma(\pt,\CC_{M})$.

We  deal with a Dirac equation over $M$ 
as an equation over another preHilbert space
$\mathcal H=(\Gamma_c(M,\Cliff_{M}^*)
\times \Gamma_c(M,\Cliff_{M}),
\langle,\rangle, \varphi)$.
Here $\LA,\RA_M$ is
the L$^2$-type pairing, for
$(\overline \psi_{M, 1},\psi_{M, 2})
\in \Gamma_c(M,\Cliff_{M}^*)\times \Gamma_c(M,\Cliff_{M})$,
\begin{equation}
	\langle\overline \psi_{M, 1},\psi_{M, 2}
        \rangle_{M}
       = \int_{M} d \mathrm{vol}_M \
\LA \overline\psi_{M, 1}, \psi_{M, 2} \RA_{\Cliff_M} 
\label{eq:pairE4}
\end{equation}
Here in $T_S$, the measure of $M$ is decomposed to
\begin{gather}
        d \mathrm{vol}_M = \rho (\det g_S)^{1/2} d^k s d^{n-k} q,
      \label{eq:detg}
\end{gather}
$d^k s = \wedge_{\alpha=1}^s d s_\alpha$, and
$d^{n-k} q = \wedge_{\dalpha=k+1}^n d q_\dalpha$.
Further in this article, we  express the preHilbert space
using the triplet with the inner product
$(\circ,\cdot)_M:=\langle\varphi \circ,\cdot\rangle_M$.
For an operator $P$ over $\Cliff_{M}$, let $Ad(P)$ be defined by
the relation,  $\langle\overline\psi_1,P \psi_2\rangle_M$
$=\langle\overline\psi_1 Ad(P), \psi_2\rangle_M$ if exists.
Further for $\psi\in\Gamma_c(M,\Cliff_{M})$, $P^*$ is defined by
$P^*\psi = \varphi^{-1}(\varphi(\psi)Ad(P))$.

Let the sheaf of the Clifford ring over $M$
 be denoted by $\CLIFF_{M}$.
As a model of (\ref{eq:gammaRm}) let
$\gamma_{M}$ be a morphism from $\Omega_M$ to $\CLIFF_M^1$.

The Dirac operator is a morphism between
the Clifford module
$$
	\Not D_M : \Cliff_{M} \to \Cliff_{M}
$$
but as a differential operator, we could extend its domain and region to,
$$
	\Not D_M : \CC_{M}^{2^{[n/2]}} \to \CC_{M}^{2^{[n/2]}}.
$$
Since $\Cliff_{M}(U)$ contains zero section,
we may consider that $\Ker ({\Not D_M})$ as a subset of
germs of $\CC_{M}^{2^{[n/2]}}$
means a subset of germs of $\Cliff_{M}$.

Then there are a set of germs
 $\{c_{M}^{a}\}_{a=1, \cdots, 2^{[n/2]}}$
of  $\Cliff_{M}(M)$ and
$ \overline{c_{M}^{a}}:=\varphi(c_{M}^{a})$
which hold relations at each point,
\begin{equation}
\LA \overline c_{M}^{a}  c_{M}^{b} \RA_{\Cliff_M} 
 = \delta^{a,b},
\quad \mbox{ for } a, b = 1, \cdots, 2^{[n/2]}.
\label{eq:SOn}
\end{equation}
A germ of solutions of Dirac equation
$\Not D_M \psi =0 $ is
expressed by $\psi = \sum_a \psi_{M}^{a} c_{M}^{a}$ for
$\psi_{M}^{a} \in \Gamma(\pt, \CC_{M})$
at a point $\pt\in M$.
Lemma \ref{lemma:SOm} gives

\begin{proposition}\label{prop:02}
There is a subsheaf $\Cliff_M^{pr} \subset\Cliff_M^{pr}$
satisfying the following:
\begin{enumerate}

\item $\Cliff_M^{pr}$ is isomorphic to  $\Theta_M$ as
vector sheaves via the following $\fj_M$, i.e., there
is a morphism $\fri: \Theta_M \to \Cliff_M^{pr}$
$(\fri(u^{(n)})=\psi_{u^{(n)}})$.

\item $\fj_M$ whose model is $\fj$ in Lemma \ref{lemma:SOm}
gives an equivalence
$\fj_M=g_M( , )$, i.e., 
for $v^{(n)}, u^{(n)}\in\Gamma(pt, \Theta_M)$,
every $\psi_{u^{(n)}}\in \Gamma(pt, \Cliff_M^{pr})$ 
satisfies
\begin{equation*}
\LA\overline \psi_{u^{(n)}} \gamma_M(g_M(v^{(n)}))
\psi_{u^{(n)}} \RA_{\Cliff_M}
       = g_M(u^{(n)}, v^{(n)}). 
\end{equation*}
We call this relation $\RR\times\SO(n)${\textrm{-representation}}
 in this article.
\end{enumerate}
\end{proposition}

Due to the Proposition,
for $\Lambda^{i}_{j}\in\Gamma(pt, \SO(n) \times \RR))$,
and $v^{(n), i}:= \Lambda^{i}_{j}e^{(n), j}\in\Gamma(pt, \Theta_M))$,
there is a pair of germ $(\psi_{e^{(n), i}})_{i=1, \cdots, n}$
in $\Gamma(pt, \Cliff_M^{pr})$ of the Clifford module
 and its dual pair
$ \overline{\psi}_{e^{(n), i}}:=\varphi_\pt(\psi_{e^{(n), i}})$
which hold a relation,
\begin{equation*}
\LA
\overline \psi_{e^{(n), \ell}} \gamma_M(g_M(v^{(n), i}))
\psi_{e^{(n), \ell}} \RA_{\Cliff_M}
       =  \Lambda^i_{\ell} \quad\mbox{(not summed over $\ell$)}.
\end{equation*}

\bigskip

Every sheaf $A_{T_S}$ over $T_S$ is determined by 
$A_{T_S} = r_{T_S, M} A_M$
for every sheaf $A_M$ over $M$ and in our conditions
these properties preserves over $T_S$.

\begin{remark}\rm{
Using a  $\CC$-valued smooth compact function
$b \in \Gamma_c(M, \CC_{M})$
over $M$ such that $b\equiv 1$ at $U \subset M$
and its support is in $M$,
$b\psi_{M}^{a}$, $b\psi_{e^{(n), k}}$ and their partners
belong to $\Gamma_c(M,\Cliff_{M})$ and $\Gamma_c(M,\Cliff_{M}^*)$.
Hereafter we assume that
$\psi_{M}^{a}$, $\psi_{e^{(n), k}}$  and their partners
are sections of $\Gamma_c(M,\Cliff_{M})$ and
$\Gamma_c(M,\Cliff_{M}^*)$
in the sense.
}
\end{remark}

The Dirac operator restricted over $T_S$ is explicitly given by
\begin{equation}
 \Not D_{T_S} =
\gamma_{T_S}(d u^\mu)(\partial_\mu +\omega_{{T_S}, \mu}).
\label{eq:DiracEqTS}
\end{equation}
where $\partial_\mu := \partial/ \partial u^\mu$
and $\omega_{{T_S},\mu}$ is a spin connection.

\section{ Submanifold Dirac Operator over $S$ in $M$}

In this section, we will define the submanifold Dirac
operator over  $S$ in $M$ and investigate
its properties.

Since $T_S$ is diffeomorphic to $\RR^n$,
$\CC_{T_S}$  is soft (Theorem 3.1 in \cite{I}).
Hence we have the following proposition.
\begin{proposition}\label{prop:41}
$\Cliff_{T_S}$ and $\CC^{2^{[n/2]}}_{T_S}$ are  soft.
\end{proposition}

\begin{proof}
$\Cliff_{T_S}$ is considered as
 a sheaf of $\CC$-vector bundle with $2^{[n/2]}$ rank.
 From the proof of Theorem 3.2 in \cite{I},
 it is justified.
\end{proof}

Due to the Proposition \ref{prop:41},
at each point $\pt$ in $S$
and for
a germ $\psi_{\pt}\in\Gamma(\pt, \Cliff_{T_S})$,
there exists $\psi_c\in\Gamma_c(T_S, \Cliff_{T_S})$
and $\psi_o\in\Gamma(T_S, \Cliff_{M})$ such that
$\psi_{\pt} = \psi_c$ and $\psi_{\pt} =\psi_o$
around $pt$.
Thus an element of
$\Gamma(\pt, \Cliff_{T_S})$ need not be distinguished
which it comes from $\Gamma_c(T_S, \Cliff_{T_S})$
or $\Gamma(T_S, \Cliff_{M})$.
From here,  every $A_{T_S}$ is identified with $A_{M}$
again.

The action of $H$ along the fiber direction,
we will continue to consider it in the framework of the
unitary representation of Clifford module and
 we wish to consider kernel of 
$\partial_\dalpha := \partial/\partial q^\dalpha$
$(\dalpha = k+1, \cdots, n)$ [\cite{Ma5} and references therein].
However $p_\dalpha:=\sqrt{-1}\partial_\dalpha$ is not
self-adjoint, $p_\dalpha^*\neq p_\dalpha$ in general
due to the existences $\rho$ in (\ref{eq:pairE4})
and (\ref{eq:detg}).

Let us follow the techniques in the pseudo-regular representation.
We introduce another preHilbert space
$\mathcal H'
\equiv ( \Gamma_c(T_S,\overwave\Cliff_{T_S}^*)
\times  \Gamma_c(T_S,\overwave\Cliff_{T_S}),
\langle,\rangle_\sa, \tilde\varphi)$ so that
$p_\dalpha$'s become  self-adjoint operators
there.
Using the half-density
(Theorem 18.1.34 in \cite{Ho}),
we construct {\textit {self-adjointization}}:
$\eta_\sa : \mathcal H \to \mathcal H'$
by,
\begin{equation*}
	\eta_\sa(\overline \psi) := \rho^{1/4} \overline\psi, \quad
	\eta_\sa(\psi) := \rho^{1/4} \psi, \quad
 	\eta_\sa(P) :=\rho^{1/4} P \rho^{-1/4}.
\end{equation*}
Here since $\rho$ does not vanish in $T_S$, $\eta_\sa$ gives
an isomorphism
$\eta_\sa:
\Cliff_{T_S}^* \times \Cliff_{T_S}
\to \overwave \Cliff_{T_S}^* \times\overwave \Cliff_{T_S}$.
Here this transformation is also essentially the same as that
in the radical Laplace operator, e.g., in Theorem 3.7 of 
\cite{He}\footnote{Our $\rho^{1/2}$ corresponds to $\delta$ in
p.261 in \cite{He}.}.
For $(\overline\psi_1,\psi_2)\in  \Gamma_c(T_S,\tilde\Cliff_{T_S}^*)
 \times  \Gamma_c(T_S,\tilde\Cliff_{T_S})$,
by letting $\tilde \varphi:= \eta_\sa \varphi \eta_\sa^{-1}$,
the pairing is defined by
\begin{equation}
	\langle\overline\psi_1,\psi_2\rangle_\sa := \int_{T_S}
 (\det g_{S})^{1/2} d^k s d^{n-k} q\
\LA  \overline\psi_1, \psi_2\RA_{\Cliff_M}.
\label{eq:pairS}
\end{equation}
Here we have the properties of $\eta_{\mathrm{sa}}$
that 1)
$\langle\circ,\cdot\rangle_\sa
=\langle\eta_\sa^{-1}\circ,\eta_\sa^{-1}\cdot\rangle_M$,
 2) for an operator $P$ of $\Cliff_{T_S}$,
$\eta_\sa(P)=\eta_\sa P \eta_\sa^{-1}$,
and 3)
 $p_\dalpha$'s themselves become self-adjoint
in $\mathcal H'$, {\textit i.e.},
$p_\dalpha =p_\dalpha^*$.

Noting $\rho=1$ at a point in $S$,
for $(\overline\psi, \psi) \in \Gamma(S, \Cliff_{T_S}^*) \times
\Gamma(S, \Cliff_{T_S})$,
we have
$$
r_{S, M}\eta_\sa(\overline\psi) = r_{S, M}\overline\psi,
\quad\mbox{and} \quad
r_{S, M} \eta_\sa(\psi) = r_{S, M} \psi.
$$

Further we have the following proposition.

\begin{proposition}\label{prop:5-1}
By letting $p_q:=a^\dalpha p_\dalpha $ for real
generic numbers $a_\dalpha$,
the projection,
\begin{equation*}
	\pi_{p_q}: \overwave \Cliff_{T_S}^* \times
             \overwave \Cliff_{T_S} \to
          \Ker( Ad(p_q))\times \Ker (p_q),
\end{equation*}
induces the projection in the preHilbert space,
 {\textit i.e.},
\begin{enumerate}
\item
For an open set $U \subset T_S$,
$\tilde\varphi|_{\Ker(p_q)}: \Gamma(U, \Ker (p_q)) \to
\Gamma(U, \Ker (Ad(p_q)))$ is isomorphic as vector space.
We simply express $\tilde\varphi|_{\Ker(p_q)}$
by $\tilde \varphi$ hereafter.

\item $\mathcal H_{p_q}:=(
\Gamma_c(T_S, \Ker (Ad(p_q)))\times
\Gamma_c(T_S, \Ker (p_q)),
\langle,\rangle_\sa,\tilde\varphi)$ is a preHilbert space.

\item  $\varpi_{p_q}:=\pi_{p_q}|_{\overwave \Cliff_{T_S}}$
induces a natural restriction of pointwise
multiplication for a point in $T_s$,
$\mathcal H^{\pt}_{p_q}:=(
\Gamma(\pt, \Ker (Ad(p_q)))\times
\Gamma(\pt, \Ker (p_q)), \cdot,
\tilde\varphi_{\pt})$ becomes a preHilbert space.
The hermite conjugate map $\tilde\varphi_{\pt}$ is
still an isomorphism.
\end{enumerate}
\end{proposition}

\begin{proof}
By letting
$\varpi_{p_q}:=\pi_{p_q}|_{\overwave \Cliff_{T_S}}$,
we have $\varpi_{p_q}=\varpi_{p_q}^2=\varpi_{p_q}^*$ in
$\mathcal H_{p_q}$.
In fact since $p_\dalpha$ is self-adjoint, $\Ker (p_q)=\Ker (p_q^*)$
and $\Ker(p_q)$ is isomorphic to
$ \Ker (Ad(p_q))$, {\textit i.e.},
$\varphi(\varpi_{p_q}\psi) = \varphi(\psi)Ad(\varpi_{p_q})$.
$\varpi_{p_q}^*\psi = \varphi^{-1}(\varphi(\psi)Ad(\varpi_{p_q}))$
gives $\varpi_{p_q}=\varpi_{p_q}^*$.
\end{proof}

\begin{remark} \label{rk:varpi}
{\rm{
We shall remark that deformation of preHilbert space by
the action of $\eta_\sa$
makes $\varpi_{p_q}$ a projection operator in the sense of
$*$-algebra. This is the essential of the scheme of the 
submanifold quantum mechanics \cite{Ma5}, which provides
non-trivial quantum mechanics \cite{MBK, GWL, MEE}.
It is absolutely non-trivial fact but the same idea appeared
in computation of Hydrogen atom in \cite{D}.
}}
\end{remark}
\bigskip
Further we consider $p_q$ as a morphism between $\CC^{2^{[n/2]}}_{T_S} 
\to \CC^{2^{[n/2]}}_{T_S}$ and its kernel 
$\Ker^{\CC} p_q \subset \CC^{2^{[n/2]}}_{T_S}$.
We are concerned with 
$r_{S, M}\Ker^{\CC} p_q \subset r_{S,M} \CC^{2^{[n/2]}}_{T_S}$,
but it is obvious that
$r_{S,M}\Ker^{\CC} p_q$ can be identified with $\CC^{2^{[n/2]}}_{S}$,
because its element is a function only of $S$.
Then we have similar relation of $\Ker^\CC p_q$ in 
Proposition \ref{prop:5-1}.

\bigskip

After we suppress a normal translation freedom in $\mathcal H_{p_q}$,
we might choose a position $q$ and make $q$ vanish.
Thus we will give our definition of the submanifold
Dirac operator.

\begin{definition}\label{SubD}{\textrm
We define the submanifold Dirac operator over  $S$ in $M$ by,
\begin{equation*}
	\Not D_{S \hookrightarrow M}
         :=r_{S, M} (\eta_\sa(\Not D_M)|_{\Ker( p_q)}),
\end{equation*}
as an endomorphism of
Clifford submodule $r_{S, M}\Ker(p_q)$ $\subset$
$r_{S, M}\overwave\Cliff_{M}$, {\it i.e.},
$$
{{\Not D_{S \hookrightarrow M}}}:
r_{S, M}\Ker(p_q) \to  r_{S, M}\Ker(p_q).
$$
Further we extend its domain and region to
$r_{S, M}\Ker^{\CC} p_q$ or $\CC^{2^{[n/2]}}_{S}$;
$$
{{\Not D_{S \hookrightarrow M}}}:
\CC^{2^{[n/2]}}_{S} \to
\CC^{2^{[n/2]}}_{S}.
$$
}
\end{definition}

Here we note that the first restriction
$|_{\Ker( p_q)}$ is as
an operator but the second one
$r_{S, M}$ is associated with a sheaf theory [\cite{I}].

\vskip 1.0  cm

In order to find the extension for
${{\Not D_{S \hookrightarrow M}}}$ over $\CC^{2^{[n/2]}}_{S}$
we need an explicit representation of the Dirac operator.
For the case that $M$ is the euclidean space, we find a natural
frame to represents the Clifford objects explicitly.
However local parameter of $M$ is not privileged in general.
Thus we introduce another Clifford ring sheaf isomorphic to
$r_{S,M}\CLIFF_{M}$ and find its explicit isomorphism
using an element of Clifford group.

Let us introduce a vector sheaf $\RR^n_S$ related to $G$-action
and a sheaf morphism $\iotaRS : \Theta_S \to \RR^n_S$
and an isomorphism $\iotaRM : r_{S, M} \Theta_M \to \RR^n_S$.
Using this, we will investigate the Clifford objects over $S$
and ones over $M$ with $r_{S, M}$ before we deal with the
Dirac operator.

Using the vector sheaf $\RR^n_S$, 
we construct a Clifford ring sheaf 
$\CLIFF(\RR^n_S)$ over $S$ generated by a linear sheaf morphism
$\gamma_{\RR^n_S}: \RR^n_S \to \CLIFF^1(\RR^n_S)$.
Similarly we could define its representation module
$\Cliff(\RR^n_S)$ and its Clifford groups $\CG(\RR^n_S)$.

We have an isomorphism
  $\iotaRM^\flat : r_{S, M}\CLIFF_M \to\CLIFF(\RR^n_S)$
and one between the Clifford groups $\CG(\RR^n_S)$ and $r_{S,M}\CG_M$.
By identifying $\CG(\RR^n_S)$ with $r_{S,M}\CG_M$,
  $\iotaRM^\flat$ is realized as $\iotaRM^\flat (c) = \tau^{-1} c \tau$
for $c \in \CLIFF_M$ and $\tau \in r_{S,M}\CG_M$.
Then we also have its representation $\Cliff(\RR^n_S)$,
 and  an isomorphism $\iotaRM^\sharp : r_{S,M}\Cliff_M \to \Cliff(\RR^n_S)$.

The $\iotaRS$ induces a ring homomorphism
$\iotaRS^{\flat} : \CLIFF_{S} \to \CLIFF(\RR^n_S)$
by its generator corresponding to  $u^{(k)}\in \Theta_S$
by $\gamma_S(u^{(k)}) \mapsto \gamma(\iotaRS(u^{(k)}))$.
Similarly we have $\iotaS^{\flat} : \CLIFF_{S} \to r_{S,M}\CLIFF_M$.
The $\iotaRS^{\flat} = \iotaRM^{\flat} \iotaS^\flat$ and $\iotaS^{\flat}$
 induce the induced and restrict representations modeling ones in \S 2
such that 
$$
{\Ind^\iotaRS}_S^{\RR^n} : \Cliff_S \to \Cliff(\RR_S^n),
\quad
{\Res^\iotaRS}_S^{\RR^n} : \Cliff(\RR_S^n) \to \Cliff_S,
$$
$$
{\Ind^\iotaS}_S^{M} : \Cliff_S \to \Cliff_M,
\quad
{\Res^\iotaS}_S^{M} : \Cliff_M \to \Cliff_S,
$$
are connected by natural relations,
$$
{\Ind^\iotaS}_S^{M} = \tau^{-1} {\Ind^\iotaRS}_S^{\RR^n}, \quad
{\Res^\iotaS}_S^{M} = {\Res^\iotaRS}_S^{\RR^n} \tau.
$$

For every $u\in \Gamma(pt, r_{S,M}\Theta_M)$,
$v \in \Gamma(pt, \Theta_S)$,
$\psi_u := \fri_M(u)$,
and  $\psi_{S,  u} := {\Res^{\iotaS}}^{M}_S \psi_{u}$,
as we showed in (\ref{eq:EK}), the Frobenius reciprocity shows
\begin{gather}
\begin{split}
 g_M(\iotaS(v), u) & =
 \LA \overline{\psi}_{u},\ \gamma_M(g_M(\iotaS_*(v)))\psi_{u} \RA_{\Cliff_M}\\
& = \LA \overline{\psi}_{S, u},\ 
        \gamma_{S}(g_{S}(v))  \psi_{S, u}
\RA_{\Cliff_S}. 
\label{eq:uvSM}
\end{split}
\end{gather}
 As in (\ref{eq:ek}), by letting 
$\psi_{\tau u} := \tau \psi_{u}\in \Gamma(pt, r_{S,M}\Cliff_{\RR^n_S})$,
we have $\psi_{S,  u} = {\Res^{\iota_{\RR^n, S}}}^{\RR^n}_S \psi_{\tau u}$
and (\ref{eq:uvSM}) becomes
\begin{gather}
\begin{split}
\LA \overline{\psi}_{\tau u},\ 
\gamma_{\RR^n}(g_{\RR^n}(\iota_{\RR^n S *}(v)))  \psi_{\tau u}
\RA_{\Cliff(\RR^n_S)} &=
\LA \overline{\psi}_{S,  u},\ \gamma_{S}(g_S(v))
 \psi_{S, u} \RA_{\Cliff_S} 
\\
& = g_S(v, \pi_{S M}(u)),
\label{eq:uvRS}
\end{split}
\end{gather}
where $\pi_{S, M}: r_{S, M}\Theta_M \to \Theta_S$ is given by
$g_M(\iotaS(v), u) = g_S(v, \pi_{S M}(u))$,
which is the simplest Frobenius reciprocity; we use its
lift to the Clifford modules.
These give the data of $\Gr^{(n, k)}_S$ and immersion $\iotaS$,
which are our purpose.

As we find relations among the Clifford objects over $S$
and ones over $M$ with $r_{S, M}$, we step to the consideration
of the Dirac operator.
In order to obtain the relation (\ref{eq:uvRS}), we will use the
Dirac operator $\Not D_{S \hookrightarrow M}$.
However we did not give its explicit representation yet.
In order to determine an explicit representation of the Dirac operator,
using $\tau \in \CG(\RR_S^n)$ which 
connects $\CLIFF(\RR^n_S)$ and $r_{S, M}\CLIFF_M$ 
as mentioned above,
we will define the Dirac operator defined over $\Cliff(\RR_S^n)$
$$
\Not D_{S \hookrightarrow M}^{\iotaRS} :=
\tau \Not D_{S \hookrightarrow M}\tau^{-1}.
$$

\begin{proposition}\label{prop:subD}
The submanifold Dirac operator of $S$ in
$M$ can be expressed by
\begin{equation}
	\Not D_{S \hookrightarrow M}^{\iotaRS}
         =
	\iotaRS^{\sharp}({\Not D}_S)
 +  \frac{1}{2}\gamma^\dalpha\iotaRM_*\Gamma_{\dalpha}.
\label{eq:DiracSE4}
\end{equation}
where $\Not D_S$ is the proper Dirac over $S$,
 $\Gamma_{\dalpha}/k$ is the mean curvature vector of $b^{(n) \dalpha}$
{\rm{[\cite{W} p.119]}} of $S$ and 
$\gamma^\dalpha :=\gamma_{\RR^n_S}(\iotaRM_*(d q^\dalpha))$.
\end{proposition}

\begin{proof}
First we note that $\eta_\sa(\Not D_{M})$ has
a decomposition,
\begin{equation*}
\eta_\sa(\Not D_{M})=
\Not{\mathbb D}_{M}^\parallel+ \Not {\mathbb D}_{M}^\perp,
\end{equation*}
where 
$\Not {\mathbb D}_{M}^\perp:=\gamma_M(d q^\dalpha)
 \partial/\partial q^\dalpha$
and $\Not {\mathbb D}_{M}^\parallel$ does not include the normal
derivative $p_\dalpha$. ${\Not D}_{M}^\perp$ vanishes
at $\Ker(p_q)$ and
at $\Ker^\CC(p_q)$.
 Due to the constructions,
$\iotaS(\gamma_S(e^{(k), \alpha}))$ and $\gamma_M(d q^\dalpha)$ 
become generator of the $\CLIFF_{T_S}$ at sufficiently vicinity of $S$.
A direct computation shows that the following relation holds
$$
r_{S, M}\left(\Not{\mathbb D}_{M}^\parallel\right) 
-\tau^{-1} \iotaRS^{\sharp}({\Not D}_S) \tau
= \frac{1}{2} 
r_{S, M}\left(\gamma_M(d q^\alpha) \Gamma_\alpha \right).
$$
The geometrical independence due to (\ref{eq:gTS}) and
direct computations give above the result.
Using $\iotaRS$ and $\iotaRM$, we have the result.
\end{proof}

\begin{remark}\label{rk:iD}
{\rm{
\begin{enumerate}
\item $\ii\iotaRS^{\sharp}({\Not D}_S)$ is a formal self-adjoint
for a $L^2$-type integral of the Clifford module over $S$ because
from the definition, $\ii{\Not D}_S$ is self-adjoint for
the integral over $S$ and $\iotaRS^{\sharp}$ is $*$-morphism.
On the other hand, $\ii\Not D_{S \hookrightarrow M}^{\iotaRS}$ is not
self-adjoint because of the extra term and the self-adjointness of
$\ii\iotaRS^{\sharp}({\Not D}_S)$.

\item Here we comment on the submanifold Dirac operators defined by
C. B\"ar [\cite{B}, Lemma 2.1]. The Dirac operator 
$\tilde D$ in \cite{B} corresponds to our
$\ii\iotaRS^{\sharp}({\Not D}_S)$ whereas
the Dirac operator $\hat D$ in \cite{B} corresponds to 
our $\ii\Not D_{S \hookrightarrow M}^{\iotaRS}$.
In \cite{Fr} the generalized Weierstrass relation is studied using
the Dirac operator which is the same as B\"ar's. Further we note that
in \cite{B}, $\tilde D$ is mainly investigated, whereas
we consider $\ii\Not D_{S \hookrightarrow M}^{\iotaRS}$,
which is not self-adjoint.

\item
Ginoux and Morel \cite{GM}, and Oussama and Zhang \cite{OZ}
dealt with the same operator $\ii\Not D_{S \hookrightarrow M}^{\iotaRS}$
but their studies started from the definition of
 $\ii\Not D_{S \hookrightarrow M}^{\iotaRS}$.
They did not mention answer why they employ the definition
in detail, at least, from viewpoint of the submanifold
quantum mechanics.

\item It is clear why 
$\ii\Not D_{S \hookrightarrow M}^{\iotaRS}$ has extra non-trivial
term. It appears due to the requirement that
the projection $\varpi_{p_q}$ should be the self-adjointness,
which is the same as the requirement  that the isomorphism
$\varphi$ should preserve for the action of $\varpi_{p_q}$.
These are essential to submanifold quantum mechanics \cite{Ma5}.
\end{enumerate}
}}
\end{remark}

Now we will give our main  theorem:
\begin{theorem}\label{th:01}
Fix the data of $\Cliff_M$ i.e.,
 its base $\LA c_{M}^{a}\RA_{a=1, \cdots, 2^{[n/2]}}$, and
a morphism $\fri: \Theta_M \to \Cliff_M^{pr}$.
Let  a point $\pt$ be in $S$ immersed in $M$.
Let $\CC^{2^{[n/2]}}_{S}$ be a sheaf of complex vector bundle over
$S$ with rank $2^{[n/2]}$.
A set of germs of $\Gamma(\pt, {\CC^{2^{[n/2]}}_{S}})$
satisfying the
submanifold Dirac equation,
\begin{equation*}
	\ii\Not D_{S \hookrightarrow M}^\iotaRS \psi
= 0 \quad \text{at}\quad \pt,
\end{equation*}
is given by $\{b_a\psi^{a} \ |\ a=1, \cdots, 2^{[n/2]},
\ b_a \in \CC \}$
 such that elements satisfy the orthonormal relation as
$\CC$-vector space;
$$
	\varphi_{\pt}(\psi^{a})\psi^{b}= \delta_{a,b}
\quad \text{at}\quad \pt.
$$
Then followings hold:

\begin{enumerate}

\item 
$\LA\psi^{a}\RA_{ a=1, \cdots, 2^{[n/2]}}$ is a base of
$\Gamma(pt, \Cliff(\RR^n_S))$.
There exists an isomorphism 
$\iotaRM^\sharp :r_{S,M}\Cliff_M \to\Cliff(\RR^n_S)$
related to $\tau \in\Gamma(\pt,r_{S,M}\CG_{M})$ satisfying
$\psi^{a}=\tau c_{M}^{a}$
$(a=1, \cdots, 2^{[n/2]})$ 
by identifying $\Cliff(\RR_S^n)$ with $r_{S,M} \Cliff_M$.
$\tau$ corresponds to an element of
$\SO(r_{S,M}\Theta_M)$ as a representative
element of $\Gr^{(n, k)}_S$.

\item 
For every $u\in \Gamma(pt, r_{S,M}\Theta_M)$, let $\psi_u := \fri_M(u) \in 
\Gamma(pt, r_{S,M}\Cliff_M^{pr})$,
$\psi_{u, S}:=\tau  \psi_{u} \in \Gamma(pt, r_{S,M}\Cliff_M^{pr})$
using $\tau$ of (1),
 and
$\overline  \psi_{u,S}$ $:=\varphi(  \psi_{u}) \tau^{-1}$
$ \in \Gamma(pt, \overline{ r_{S,M}\Cliff_M^{pr}})$.
Then for every $v\in \Gamma(pt, r_{S,M}\Theta_S)$, the following relation holds:
\begin{equation}
\LA \overline\psi_{u, S} [\iotaRS^{\flat}
(\gamma_{S}(g_S(v)))] \psi_{u, S} \RA_{\Cliff(\RR^n_S)}
           = g_M( \iotaS(v), u).
\end{equation}
This value brings us the local data of immersion $\iotaS$.
\end{enumerate}

\end{theorem}

\begin{proof}
Since $\Not D^\iotaRS_{S \hookrightarrow M}$ is the
$2^{[n/2]}$ rank first order differential operator and
has no singularity over $S$ due to the construction,
a germ of its kernel in $\Gamma(\pt, \CC_{S}^{2^{[n/2]}})$
 is given by $2^{[n/2]}$ dimensional
vector space at each point of $S$.
Since $\Not D^\iotaRS_{S \hookrightarrow M}$
is defined as an endomorphism of $\Ker^\CC(p_q) \approx
\CC^{2^{[n/2]}}_S$.
The kernel of the Dirac operator,
$\Ker^\CC({\Not D}_{S \hookrightarrow M}^\iotaRS)$
of $\CC^{2^{[n/2]}}_S$ has an injection into
 $\Cliff(\RR^n_S)$.
There exist $\tau \in \CG(\RR^n_S)$ such that
 ${\iotaRM^\flat}^{-1} :
 r_{S, M}\Ker^\CC({\Not D}_{S \hookrightarrow M}^\iotaRS)$
$\to \Cliff_M$.

Let $\Not{\mathbb D_S}^\perp := \tau^{-1} \gamma^\dalpha \partial_\dalpha \tau$
at $S$.
From the construction, we have
$$
	\iotaS_* \Not D_{S \hookrightarrow M} + \Not{\mathbb D_S}^\perp
       = r_{S, M}(\eta_\sa(\Not D_{D})) .
$$
Hence $\Ker^\CC(\Not D^\iotaRS_{S \hookrightarrow M})$ is
a subset of 
 a kernel of $\tau (r_{S, M}(\eta_\sa({\Not D}_{M})) )\tau^{-1}
\subset r_{S, M} \CC_{M}^{2^{[n/2]}}$.

Noting Proposition \ref{prop:5-1}, $\tilde \varphi_{\pt}$ is an isomorphism
and $\mathcal H^{\pt}_{p_q}$ gives
(5.2) and (5.3). Thus we prove them.
\end{proof}

\begin{remark}
{\rm{
\begin{enumerate}
\item 
The finial result does not depend upon a choice of $\iotaRS$.

\item
This theorem is based upon
the Frobenius reciprocity of Clifford ring sheaves
 on category of differential geometry
as shown in (\ref{eq:uvRS}) and (\ref{eq:uvSM}).
We have
compared ${\Ind^\iotaRS}_S^{\RR^n}  \Cliff_S$, which is obtained by
using the Dirac operator, with $\Cliff_M$ as each germ in Theorem 5.1.

\item We have assumed that $M$ and $S$ are homeomorphic to
$\RR^n$ and $\RR^k$ respectively.
However as our arguments are local, the theorem could be
extended to spin manifolds $S$ and $M$ under assumptions
on the group action if there is no geometrical obstruction.

\item With Remark \ref{rk:varpi} and \ref{rk:iD} (4),
 it is obvious that the submanifold
Dirac operator given in  submanifold quantum mechanics
represents local immersed geometry. 
Its essential is that the restriction of the 
Dirac operator preserving $\varphi$ in Definition 5.1
consists with the Frobenius
reciprocity. It is the answer of the question mentioned in
Introduction.

\item If $M$ and $S$ have natural parameterization 
$(x^{i})_{i=1, \cdots, n}$ and $(s^\alpha)_{\alpha=1, \cdots, k}$
and $S$ is  an analytic submanifold such that 
$$
x^{ i}(s) = \int^s_S d x^{ i}(s)
$$
represents an immersion $S$ in $M$, it can be expressed as
$$
x^{i}(s) = \int^s_S g_{S, \alpha,\beta}
\LA\overline\psi_{\partial_{ x^{i}}, S} [\iotaS^{\sharp}
(\gamma_{S}(ds^\alpha))] 
\psi_{\partial_{ x^{i}}, S}\RA_{\Cliff(\RR^n_S)} d s^\beta,
$$
where $\partial_{ x^{i}}:= \partial /\partial x^{i}$
using above $\psi$. This is the generalized Weierstrass relation.

\item When $M\equiv \RR^n$ and  $k=2$, the theorem is reduced to
the generalized Weierstrass relation \cite{Fr, Ko, KL, PP}.
In the case, $\iotaRS$ is properly determined and identify $\RR_S^n$
with $\RR^n$.
These are closely related to the two-dimensional integrable system.
Especially, when $\Not D^\iotaRS_{S \hookrightarrow M}$ is identified with
${\Not D}_S$ and $\bar \partial$, which correspond to
minimal surface cases,
it becomes original Weierstrass relation [\cite{E}, p.260-7].

\item
As mentioned in \cite{Ma4}, we can put the Frenet-Serret torsion 
field into the Dirac operator.

\item For $k=1$ case, Theorem is mere the Frenet-Serret relation
\cite{Ma1, Ma2}.

\item As we showed in \cite{Ma1, Ma2}, the Dirac operator also
might give the global properties of the immersion of $S$, 
i.e., its topological properties,
though we mentioned only local properties in this article.
Thus we should investigate the global properties using the 
submanifold Dirac operator as generalization of \cite{Ma1, Ma2}
in future.

\item When $S$ is a conformal surface,
we may consider the relations along the line of arguments of
\cite{BF, Ko, Ta1, Ta2, Ta3, PP}. For example, we could classify 
the immersions using the Dirac operator.
Furthermore  when $S$ has holomorphic properties,
we also may give similar arguments.

\end{enumerate}
}}
\end{remark}

\bigskip

{\bf{Acknowledgment}}

The author thanks Professor K. Tamano, Professor N. Konno and 
Dr. H. Mitsuhashi for encouragements on this work, especially
Dr. H. Mitsuhashi for his lecture on Frobenius reciprocity.


\begin{thebibliography}{ABS}

\bibitem{ABS}
\by{M. F. Atiyah, R. Bott and A. Shapiro}
    \paper{Clifford modues}
\jour{Topology}
\vol{3} \yr{1964}, \pages 3-38.


\bibitem{B}
\by{C. B\"ar}
 \paper{Extrinsic Bounds for Eigenvalues of the Dirac operator}
\jour{Ann. Glob. Anal. Geom.} \vol{16} \yr{1998} \pages{573--596}.

\bibitem{BF}
\by{L. V. Bogdanov and E. V. Ferapontov}
\paper{Projective differential geometry of higher
reductions of the two-dimensional Dirac equation}
\jour{J. Geom. Phys.} \vol{52} \yr{2004} \pages{328-352}.

\bibitem{BGV}
\by{N. Berline, E. Getzler, M. Vergne}
    \book{Heat kernels and Dirac operators}
Springer, Berlin, 1996.

\bibitem{C}
\by{C. Chevalley}
\book{The algebraic theory of spinors and Clifford algebras}
Springer, Berlin, 1997.

\bibitem{D}
\by{P. A. M. Dirac}
\book{The principles of Quantum Mechanics} fourth edition,
Oxford Univ. Press, Oxford, 1958.

\bibitem{E}
\by{K. P. Eisenhart}
\book{A treatise on the differential geometry of curves and surfaces}
Ginn and Company, Boston, 1909.

\bibitem{En}
 M. Encinosa,
\paper{Electron wave functions on $T^2$ in a static magnetic field
of arbitrary direction}
\jour{Physica E: Low-dimensional Systems and Nanostructres}
    \vol{28} \yr{2005} 209--218.


\bibitem{Fr}
\by{T. Friedrich}
 \paper{On the spinor representation of surfaces in Euclidean 3-space}
\jour{J. Geom. Phys.} \vol{28} \yr{1998} \pages{143--157}.

\bibitem{GM}
N. Ginoux and B. Morel,
\paper{On the eigenvalue estimates for the submanifold Dirac operator}
\jour{Int. J. Math.} \vol{13} \yr{2002} 533-548.

\bibitem{GW}
\by{R. Goodman and N. R. Wallach}
\book{Representations and Invariants of the Classical Groups}
Cambridge Univ. Press, Cambridge, 2003


\bibitem{GWL}
\by{J. Gravesen, M. Willatzen, and L. C. Lew Yan Voon}
\paper{Schr\"odinger problems for surfaces of revolution--the 
finite cylinder as a test example}
\jour{J. Math. Phys.} \vol{46} \yr{2005} 012107.

\bibitem{Ho}
\by{L. H\"ormander}
\book{The analysis of linear partial differential
operators III}
Springer-Verlag Berlin, 1985.

\bibitem{He}
\by{S. Helgason}
\book{Differential geometry, Lie groups, and symmetric spaces}
 \publ{Academic Press} \yr{1978}

\bibitem{I}
\by{B. Iversen}
    \book{Cohomology of Sheaves}
    Springer-Verlag, 
  1986.

\bibitem{Ke}
K. Kenmotsu,
\paper{Weierstrass formula for surfaces of prescribed mean curvature}
\jour{Math. Ann.}
\vol{245} \yr{1979}, \pages{89-99}.

\bibitem{Ko}
B. G. Konopelchenko,
\paper{Weierstrass representations for surfaces in 4D spaces
and their integrable deformations via DS hierarchy}
\jour{Ann. Global Analysis and Geom.}
\vol{16} \yr{2000}, \pages{61-74}.

\bibitem{KL}
\by{B. G. Konopelchenko and G. Landolfi}
\paper{Generalized Weierstrass representation for surfaces
in multi-dimensional Riemann spaces}
\jour{J. Geom. Phys.} \vol{29} \yr{1999} \pages{319--333}.

\bibitem{Ma1}
S. Matsutani,
         \paper{Anomaly on a submanifold system -New index theorem
        related to a submanifold system}
   \jour{J. Phys. A} \vol{28} \yr{1995} 1399-1412.

\bibitem{Ma2}
S. Matsutani,
         \paper{Immersion anomaly of Dirac operator of 
                of surface in $\RR^3$}
   \jour{Rev. Math. Phys.} \vol{11} \yr{1999} 171-186.

\bibitem{Ma3}
S. Matsutani, \paper{Generalized Weierstrass relation for a submanifold $S$
in $\EE^n$ arising from the submanifold Dirac operator}
   to appear in \jour{Adv. Stud. Pure Math.}.

\bibitem{Ma4}
S. Matsutani,
         \paper{Submanifold Dirac operators with torsion}
   \jour{Balkan J. Geom. and Its Appl.} \vol{9} \yr{2004} 1-5.

\bibitem{Ma5}
S. Matsutani,
         \paper{On the essential algebraic aspect
                of submanifold quantum mechanics}
   \jour{J. Geom. and Symm. in Phys.} \vol{2} \yr{2004} 18-26.

\bibitem{MBK}
G. J. Meyer, R. H. Blick and I Knezevic,
\paper{Curvature-Dependent Conductance Resonances in Quantum
Cavities} 2005.

\bibitem{MEE}
L. Mott, M. Encinosa, and B. Etemadi,
\paper{A numerical study of the spectrum and eigenfunctions
on a tubular arc}
\jour{Physica E: Low-dimensional Systems and Nanostructres}
    \vol{25} \yr{2005} 532--529.

\bibitem{MT}
S. Matsutani and H. Tsuru,
         \paper{Physical relation between quantum mechanics and 
           solitons on a thin elastic rod}
   \jour{Phys. Rev. A} \vol{46} \yr{1992} 1144-1447.

\bibitem{OZ}
H. d'Oussama and X. Zhang,
\paper{Lower bounds for the eigenvalues of the Dirac operator,
part II. The submanifold Dirac operator}
\jour{Ann. Global Anal. Geom.} \vol{19} \yr{2001} 163-181.


\bibitem{PP}
  F. Pedit and U. Pinkall,
\paper{Quaternionic Analysis on Riemann Surfaces and Differential Geometry}
\jour{Doc. Math. J. DMV} Extra Vol. ICM II
(1999), 389-400.

\bibitem{S}
J-P. Serre,
\book{Linear Representations of Finite Group}
Springer,
1977.

\bibitem{Ta1}
 I. A. Taimanov,
\paper{The Weierstrass representation of closed surfaces in $ R\sp 3$}
\jour{Funct. Anal. Appl.} \vol{32} (1998), 258--267.

\bibitem{Ta2}
 I. A. Taimanov,
\paper{Surfaces in the four-space and the Davey-Stewartson equations}
\jour{J. Geom. Phys.} \vol{56} (2006), 1235--1256.

\bibitem{Ta3}
 I. A. Taimanov,
\paper{Two-dimensional Dirac operator and the theory of surfaces}
\jour{Russian Math. Surveys} \vol{61} (2006), 79--159.

\bibitem{V}
 V. V. Varlamov,
\paper{Generalized Weierstrass representation for surfaces in terms of 
Dirac-Hestenes spinor field}
\jour{J. Geom. Phys.} \vol{32} (2000), 241--251.

\bibitem{W}
T. J. Willmore,
\book{Riemannian geometry}
Clarendon press, Oxford
1993.
\end{thebibliography}

\bigskip

\leftline{{Shigeki Matsutani}}

\leftline{{e-mail:rxb01142@nifty.com}}

\leftline{{8-21-1 Higashi-Linkan,}}

\leftline{{Sagamihara 228-0811}}

\leftline{{JAPAN}}

\end{document}